\documentclass[a4paper,12pt]{amsart}
\usepackage{amsfonts}
\usepackage{amssymb}
\usepackage{ifthen}
\usepackage{graphicx}
\usepackage[usenames]{color}

\makeatletter
\@namedef{subjclassname@2020}{%
  \textup{2020} Mathematics Subject Classification}
\makeatother

\numberwithin{equation}{section}

\newcounter{alphabet}

\newenvironment{Thm}[1][]{\refstepcounter{alphabet}%
\bigskip%
\noindent%
{\bf Theorem \Alph{alphabet}}%
\ifthenelse{\equal{#1}{}}{}{ (#1)}%
{\bf .} \itshape}{\vskip 8pt}

\newcommand{\IR}{{\mathbb R}}

\newcommand{\IC}{{\mathbb C}}
\newcommand{\ID}{{\mathbb D}}
\newcommand{\IB}{{\mathbb B}}

\def\be{\begin{equation}}
\def\ee{\end{equation}}

\newcommand{\bee}{\begin{enumerate}}
\newcommand{\eee}{\end{enumerate}}

\newcommand{\blem}{\begin{lem}}
\newcommand{\elem}{\end{lem}}
\newcommand{\bthm}{\begin{thm}}
\newcommand{\ethm}{\end{thm}}
\newcommand{\bcor}{\begin{cor}}
\newcommand{\ecor}{\end{cor}}
\newcommand{\beg}{\begin{examp}}
\newcommand{\eeg}{\end{examp}}
\newcommand{\begs}{\begin{examples}}
\newcommand{\eegs}{\end{examples}}
\newcommand{\bdefe}{\begin{defn}}
\newcommand{\edefe}{\end{defn}}
\newcommand{\bprob}{\begin{prob}}
\newcommand{\eprob}{\end{prob}}
\newcommand{\bques}{\begin{ques}}
\newcommand{\eques}{\end{ques}}
\newcommand{\bei}{\begin{itemize}}
\newcommand{\eei}{\end{itemize}}
\newcommand{\bcon}{\begin{conj}}
\newcommand{\econ}{\end{conj}}
\newcommand{\bcons}{\begin{conjs}}
\newcommand{\econs}{\end{conjs}}
\newcommand{\bprop}{\begin{prop}}
\newcommand{\eprop}{\end{prop}}
\newcommand{\br}{\begin{rem}}
\newcommand{\er}{\end{rem}}
\newcommand{\bo}{\begin{obser}}
\newcommand{\eo}{\end{obser}}
\newcommand{\bos}{\begin{obsers}}
\newcommand{\eos}{\end{obsers}}
\newcommand{\bpf}{\begin{proof}}
\newcommand{\epf}{\end{proof}}
\newcommand{\ba}{\begin{array}}
\newcommand{\ea}{\end{array}}
\newcommand{\beq}{\begin{eqnarray}}
\newcommand{\beqq}{\begin{eqnarray*}}
\newcommand{\eeq}{\end{eqnarray}}
\newcommand{\eeqq}{\end{eqnarray*}}

\newtheorem{thm}{Theorem}[section]
\newtheorem{lem}{Lemma}[section]

\begin{document}
\title{Zygmund type spaces of harmonic functions on the real unit ball}

\author[X. Fu]{Xi Fu}
\address{Xi Fu, School of Mathematics Physics and Statistics\\
Shanghai Polytechnic University\\
Shanghai 201209, P. R. China. }
\email{fuxi@sspu.edu.cn}

\author[A. Rasila]{Antti Rasila}
\address{Antti Rasila, Department of Mathematics with Computer Science\\
Guangdong Technion Israel Institute of Technology\\
Shantou, Guangdong 515063, P. R. China\\ and
Department of Mathematics, Technion-Israel Institute of Technology, Haifa
3200003, Israel.
}
\email{antti.rasila@gtiit.edu.cn, antti.rasila@iki.fi}

\author[X. Xie]{Xiaoqiang Xie}
\address{Xiaoqiang Xie, School of Mathematics Physics and Statistics\\
Shanghai Polytechnic University\\
Shanghai 201209, P. R. China. }
\email{xqxie@sspu.edu.cn}



\date{}

\begin{abstract}
In this paper, we obtain several characterizations of harmonic Zygmund type spaces on the real unit ball $\IB$ of $\IR^n$. First, we characterize the harmonic Zygmund space $\mathcal{HZ}^\alpha$ ($0<\alpha\leq2$) in terms of Zygmund type conditions. Our result extends \cite[Theorem 3.4]{AC} to the setting of the real unit ball $\IB$. Second, we establish an analogue of the Holland-Walsh characterization of the harmonic Zygmund space  $\mathcal{HZ}$. Finally, an integral characterization of $\mathcal{HZ}^\alpha$ is also discussed. All obtained results can be viewed as counterparts of the known results for Bloch spaces.
\end{abstract}

\maketitle


\section{Introduction and  main results}
\label{}

Let $\IR^n$ be the $n$-dimensional real vector space, and let $e_1,\ldots,e_n$ denote the standard orthonormal basis in $\IR^n$. For $x=(x_1,\ldots,x_n), y=(y_1,\ldots,y_n) \in \IR^n$,   the inner product $\langle \cdot , \cdot\rangle$ of $x,y$ is defined as $\langle x,y\rangle= x_1y_1+\cdots+x_ny_n$ so that the Euclidean norm of $x$ is given by
$$|x|=\sqrt{\langle x,x\rangle}=\sqrt{x_1^2+\cdots+x_n^2}. $$
For $a\in \IR^n$ and $r>0$, let $\IB(a,r)$ and $\mathbb{S}(a,r)$ be the Euclidean ball  and the Euclidean sphere in $\IR^n$ with center $a$ and radius $r$, respectively. In particular, we use the notations $\IB_r:= \IB(0,r)$, $\IB:= \IB(0,1)$,
$\mathbb{S}:=\mathbb{S}(0,1)$,  and
$\overline{\IB}:=\IB \cup \mathbb{S} $ the closure of $\IB$.  We denote by $h(\IB)$ the class of
all harmonic functions on $\IB$ and by $h(\overline{\IB})$ the subclass of $h(\IB)$ consisting of all harmonic functions that are continuous on $\overline{\IB}$.

For $\alpha>0$, we define the {\it harmonic $\alpha$-Bloch space }  $\mathcal{HB}^\alpha$ as the space of functions $u\in h(\IB)$ such that
$$B_u^\alpha=\sup_{x\in \mathbb{B}}(1-|x|^2)^\alpha  |\nabla u(x)| <\infty,
$$  where $\nabla$ denotes the gradient operator.
The {\it harmonic $\alpha$-Zygmund space} $\mathcal{HZ}^\alpha$ consists of all $u\in h(\IB)$ for which
\begin{equation}\label{zygmund}Z_u^\alpha= \sup_{x\in \mathbb{B}}(1-|x|^2)^\alpha |\nabla^2u(x)| <\infty,
\end{equation} where $|\nabla^2u(x)|$ denotes the $\ell_1$ norm of the Hessian matrix defined by
$$|\nabla^2u(x)|= \sum^n_{i=1}\sum^n_{j=1} \Big|\frac{\partial^2u}{\partial x_i\partial x_j} (x) \Big|.$$

It is well-known that $\mathcal{HB}^\alpha$ and  $\mathcal{HZ}^\alpha$  are Banach spaces with the norm $$\|u\|_{\mathcal{HB}^\alpha}=|u(0)|+B_u^\alpha ~ \;\;\;\ \mbox{and} ~ \|u\|_{\mathcal{HZ}^\alpha}=|u(0)|+Z_u^\alpha,$$  respectively. In view of \eqref{zygmund},  $u\in \mathcal{HZ}^\alpha$ if and only if  $u_i= \frac{\partial u}{\partial x_i}\in \mathcal{HB}^\alpha$ for $i=1,2,\ldots,n$.  In particular, when $\alpha=1$, the spaces  $\mathcal{HB}^1$ and $\mathcal{HZ}^1$ are the classical harmonic Bloch space $\mathcal{HB}$  and harmonic Zygmund space $\mathcal{HZ}$, respectively.

Let $\alpha\in(0,1]$ and $\beta\in (0,2]$. A function $f$ is said to belong to the {\it $\alpha$-H\"{o}lder
class} $\Lambda^\alpha(\IB)$ ~(resp. $\Lambda^\alpha(\mathbb{S})$ ) if there exists a constant $C=C_f$ such that \begin{eqnarray}|f(x)-f(y)|\leq C|x-y|^\alpha, ~ x,y \in \IB ~(\mbox{resp.}  ~x,y \in \mathbb{S} ).\end{eqnarray}
 $f$ is said to belong to the {\it $\beta$-Zygmund class} $Z^\beta(\IB)$ if it satisfies the following Zygmund type condition:
  \begin{eqnarray} Z_f(x,y)=|f(x+y)+f(x-y)-2f(x)|\leq C'|y|^\beta,\end{eqnarray} for all  $x\in \IB$ and $y$ with $x\pm y\in \IB$ where $C'=C'_f$ is a constant.

For investigations of the $\alpha$-H\"{o}lder  class and the $\beta$-Zygmund class of  analytic (harmonic) functions on the unit disk $\ID$, the unit ball $\IB$ and the upper half space we refer to \cite{Du, Dy1, Dy2, FRX,  Kr} and the references therein. An important result of Hardy and Littlewood states that an analytic function $f$ on $\ID$ belongs to $\Lambda^\alpha(\ID)$ (resp. $Z^\beta(\ID)$) if and only if $f \in \mathcal{B}^{1-\alpha}$ (resp. $f \in \mathcal{Z}^{2-\beta}$), where $\mathcal{B}^\alpha$ and $\mathcal{Z}^\beta$ are the analytic $\alpha$-Bloch space and $\beta$-Zygmund space, respectively (cf. \cite{BD, Du, Zy}).

Suppose that $u\in h(\overline{\IB})$ with the boundary function $u^*=u|_{\mathbb{S}}$. Then  $u$ can be represented as the Poisson integral of $u^*(\xi)$ on $\mathbb{S}$ as
\begin{equation}\label{equ2.2}
u(x)=P[u^*](x)=\int_{\mathbb{S}}P(x,\xi)u^*(\xi)~d\sigma(\xi), \;\;\ x \in \IB,
\end{equation} where  $P(x,\xi)=\frac{1-|x|^2}{|x-\xi|^n}$ is the
 Poisson kernel for $\IB$ and $d\sigma$ is the normalized surface measure on $\mathbb{S}$.

 It is well known that the condition $u^*\in \Lambda^\alpha(\mathbb{S})$ guarantees  its harmonic extension $P[u^*]\in \Lambda^\alpha(\IB)$ when  $\alpha \in (0,1)$ (cf. \cite{Zy}). However, the result is not true for the case of $\alpha=1$.  Arsenovi\'{c}, Koji\'{c} and Mateljevi\'{c} \cite{AKM}
showed that $u^*\in \Lambda^1(\mathbb{S})$ implies $P[u^*]\in \Lambda^1(\IB)$ provided that $P[u^*]$ is a $K$-quasiregular mapping.

Recall that, by a result of  Dyakonov \cite{Dy1}, an analytic function $f$ on $\ID$, continuous up to the boundary $\mathbb{T}=\partial \ID$, is  $\alpha$-H\"{o}lder continuous on $\ID$ if and only if its modulus $|f|$ is $\alpha$-H\"{o}lder continuous on $\mathbb{T}$, i.e., $|f| \in \Lambda^\alpha(\mathbb{T})$,  and the following condition is satisfied:
$$ \big||f(re^{i\theta})|-|f(e^{i\theta})|\big|\leq C(1-r)^\alpha, ~~r\in (0,1), ~\theta\in [0,2\pi),$$
where $C$ is a constant. A simple proof of this result was given by Pavlovi\'{c} \cite{Pa1}. Pavlovi\'{c} \cite{Pa} has also extended Dyakonov's result to the setting of harmonic functions on $\IB$. For the further discussion on this topic, we refer to \cite{CHRW,CPR, KRD, KM, LW, O}.

\subsection{Zygmund type conditions for harmonic functions}
Recently, Aljuaid and Colonna \cite{AC} discussed an analogue of this problem in the setting of harmonic mappings on $\ID$ and obtained the following result.

\begin{Thm}
Let $f$ be a harmonic mapping on $\ID$ and continuous up to $\mathbb{T}$. Then $f\in \mathcal{HZ} $ if and only if
\begin{eqnarray*}
\sup_{t>0}\frac{|f(e^{i(\theta+t)})+ f(e^{i(\theta-t)})-2f(e^{i\theta})|}{t}<\infty.
\end{eqnarray*}
\end{Thm}

In view of the above-mentioned works, it is natural to ask if Zygmund type conditions exist to characterize $\mathcal{HZ}^\alpha$, and our first result provides such a characterization. Before stating our result, we need to introduce the concept of a {\it symmetric mapping } on $\mathbb{S}$.

For $x'\in \mathbb{S}$,  define \begin{equation}\sigma_{x'}(y')=2\langle x',y'\rangle x'-y', ~~ y'\in \mathbb{S}.\end{equation} It is easy to check that $\sigma_{x'}$ is an involution of $\mathbb{S}$ satisfying  $$|x'-y'|=|x'-\sigma_{x'}(y')| ~~\;\; \mbox{and} \;\;\ \langle x' , y'\rangle= \langle x' , \sigma_{x'}(y')\rangle.$$ We say that {\it $y'$ and $\sigma_{x'}(y')$ are  symmetric points on  $\mathbb{S}$ with respect to $x'$}.
\smallskip

\begin{thm}\label{th1}
Let $\beta\in (0,2)$ and $u\in h(\overline{\IB})$ with the boundary function $u^*=u|_{\mathbb{S}}$. Then the following conditions are equivalent:
\smallskip

\noindent $(a)$  $u\in \mathcal{HZ}^{2-\beta}$;
\smallskip

\noindent  $(b)$ There exists a constant $C$ such that $$|u^*(b')+u^*(c')-2u^*(a')|\leq C|a'-b'|^\beta, ~ a', b', c'\in \mathbb{S},  c'=\sigma_{a'}(b');$$
\smallskip
\noindent  $(c)$  There exists a constant $C'$ such that $$|u(x+y)+u(x-y)-2u(x)|\leq C'|y|^\beta, ~x\in \IB ~ \mbox{and} ~y ~with~ x\pm y\in \IB .$$
\end{thm}

\smallskip

Let $x\in \IB$ with $|x|=r$. Denote by  $T=T_x r\mathbb{S}$  the tangent hyperplane at $x$ to the sphere $r\mathbb{S}$. For the critical case $\beta=2$, we have the following result.

\begin{thm}\label{th2}
Let $u\in h(\overline{\IB})$ with the boundary function $u^*=u|_{\mathbb{S}}$. Assume that
there exists a constant $C$ such that $$|u^*(b')+u^*(c')-2u^*(a')|\leq C|a'-b'|^2, ~ a', b', c'\in \mathbb{S},  c'=\sigma_{a'}(b'),$$ Then $|\nabla^2(u|_{T})|$ is bounded.

\smallskip
 \noindent Moreover, if
\begin{equation} |\nabla^2u(x)|\leq L |\nabla^2(u|_{T})(x)|,~~x\in \IB, \end{equation}for some constant $L$, then $u\in Z^2(\IB)$.
\end{thm}

\smallskip

\noindent {\bf Remark 1.} In \cite{AKM}, the authors proved that if $u^*\in \Lambda^1(\mathbb{S})$, then $|\nabla(u|_{T})|$ is bounded. By the quasiregularity of $u$, it follows immediately that  $|\nabla u|\leq C|\nabla(u|_{T})|$ and thus $u\in \Lambda^1(\mathbb{B})$. Therefore, Theorem 1.2 can be viewed as a generalization of this result to the second iterated difference setting.

\medskip

\subsection{Holland-Walsh criteria for harmonic function spaces} For $\mu,\nu\geq 0$ and $f\in \mathcal{C}(\IB)$, if there exists a constant $C$ such that
\begin{equation}\label{wght_lptz}
(1-|x|^2)^{\mu}(1-|y|^2)^{\nu}|f(x)-f(y)|\leq C|x-y|,
\end{equation}
for all $x,y \in \IB$, then we say that $f$ satisfies a {\it weighted Lipschitz condition}  of indices $(\mu,\nu)$ (cf. \cite{RK1,RK}).

In 1986,  a celebrated criterion of  Holland and Walsh \cite{HW} states that an analytic function $f$ on $\ID$ belongs to
  the Bloch space $\mathcal{B}$ if and only if it satisfies the weighted Lipschitz condition
 with indices $(1/2,1/2)$, namely, \eqref{wght_lptz} with $\mu=\nu=1/2$. This result can be regarded as a derivative-free characterization of  Bloch
spaces. Since then, the question of obtaining  derivative-free characterizations of Bloch,
 $\alpha$-Bloch, and Besov
spaces of analytic and harmonic functions has attracted much
attention. In \cite{RT}, Ren and Tu extended the Holland-Walsh characterization
to the several variables case on the unit ball of $\IC^n$.   By adding the restriction ``$w$ belongs to a small neighborhood of $z$",  Li and Wulan \cite{LW1} characterized the $\alpha$-Bloch space in terms of a weighted Lipschitz condition
 of indices $(\beta,\alpha-\beta)$.  See e.g. ~\cite{Ma, Pa2, RK1, RK,zh} for more discussion on characterizations of Bloch, Besov and Zygmund spaces. This work also continues our previous study \cite{FRX} where problems of this type in the half-space  case were considered.

Our second result is an analogue of the Holland-Walsh characterization in the setting of  $\mathcal{HZ}$.

\begin{thm}\label{th3}
Let $u\in h(\IB)$.  Then
$u \in \mathcal{HZ}$ if and only if for all $x \in \IB$ and $y$ with $x\pm y\in \IB$ such that
\begin{eqnarray}
I=\sup_{y\neq 0}~
(1-|x+y|^2)^{1/2}(1-|x-y|^2)^{1/2}\frac{Z_u(x,y)}{|y|^2}<\infty.
\end{eqnarray}
\end{thm}

By adding a restriction $x\pm y \in E(x,r)$, we obtain the following result.

 \begin{thm}\label{th4}
Let $0< \beta\leq \alpha$, $0<r<1$, and let $u\in h(\IB)$.  Then
$u\in \mathcal{HZ}^\alpha$ if and only if
\begin{equation}J=
\sup_{x\pm y \in E(x,r), y\neq 0}~
(1-|x+y|^2)^{\beta}(1-|x-y|^2)^{\alpha-\beta}\frac{Z_u(x,y)}{|y|^2}<\infty.
\end{equation}
\end{thm}

Let $\mathbf{B}$  be the complex unit ball of $n$-dimensional complex vector space
$\IC^n$. For $a\in \mathbf{B}$, let $\phi_a$ be the canonical biholomorphic self map of $\mathbf{B}$  satisfying
$\phi_a(0)=a$ and $\phi_a(a)=0$. In \cite{LW2}, Li and Wulan discussed integral characterizations of Bloch space and obtained that an analytic function $f$ of $\mathbf{B}$ belongs to the Bloch space $\mathcal{B}$ if and only if
\begin{eqnarray*} \sup_{a\in \mathbf{B}}\int_{\mathbf{B}}|\widetilde{\nabla} f(z)|^p(1-|z|^2)^{-n-1}(1-|\phi_a(z)|^2)^{n+1}dV(z)<\infty, \end{eqnarray*}
where $\widetilde{\nabla} f$ is the invariant gradient of $f$ and $dV$ denotes the volume measure in $\mathbf{B}$.

For $a\in \IB$, denote by  $\varphi_a$ the M\"{o}bius transformation of $\IB$ with $\varphi_a(0)=a$ and $\varphi_a(a)=0$ (see details in Section 2).  Our final result establishes a counterpart of \cite[Theorem 3]{LW2} for $\mathcal{HZ}^\alpha$ in terms of $\varphi_a$ as follows.

\begin{thm}\label{th5}
Let $\alpha>0$, $s>n-1$, $0<p<\infty$ and $u\in h(\IB)$.  Then $u\in \mathcal{HZ}^\alpha$ if and only if
\begin{eqnarray} \sup_{a\in \IB}\int_{\IB}|\nabla^2u(x)|^p(1-|x|^2)^{p\alpha-n}(1-|\varphi_a(x)|^2)^sdv(x)<\infty, \end{eqnarray} where $dv$ denotes the volume measure in $\IB$.
\end{thm}

 In Section 2, some
necessary terminology and notation will be introduced.  The proofs of the main results will be presented in Sections 3 and 4.

\section{ Preliminaries}

We shall start this section with some preliminaries and auxiliary results. Throughout this paper, positive constants are denoted by $C$, and these constants may differ from one occurrence to the other. For
nonnegative quantities $X$ and $Y$, the relation $X \lesssim Y$ means that $X$ is
dominated by $Y$ times some positive constant. The notation
$X\thickapprox Y$ means that $Y\lesssim X\lesssim Y$. \smallskip

Let $x \in \IR^n$, write $x$ as $x=|x|x'$, $x' \in \mathbb{S}$. For $a,b \in \IR^n$, let
$$[a,b]=\Big||a|b-a'\Big|.$$
Then the symmetry lemma shows that
$$[a,b]=[b,a].$$

For  $a\in \IB$, the M\"{o}bius transformation $\varphi_a$ from $\IB$
onto $\IB$ is given by \begin{eqnarray*}\varphi_a(x)=\frac{|x-a|^2a-(1-|a|^2)(x-a)}{[x,a]^2}, \; x\in \IB.\end{eqnarray*}
It is known that $\varphi_a$  is an involution of $\IB$ such that $\varphi_a(0)=a$ and $\varphi_a(a)=0$.

By straightforward calculations, one can obtain $$|\varphi_a(x)|=\frac{|x-a|}{[x,a]}$$ and
$$1-|\varphi_a(x)|^2=\frac{(1-|x|^2)(1-|a|^2)}{[x,a]^2}. $$
For basic properties of M\"{o}bius transformation $\varphi_a$, see e.g. \cite{Ah, RK1}.

\smallskip
For $a\in \IB$ and $r\in (0,1)$,   the {\it pseudo-hyperbolic ball}
with center $a$ and radius $r$ is defined as  $$E(a,r)=\{x\in
\mathbb{B}: |\varphi_a(x)|<r\}.$$
Indeed, $E(a, r)$ is a Euclidean ball with center $c_a$ and radius
$r_a$ given by
\begin{eqnarray*}c_a=\frac{(1-r^2)a}{1-|a|^2r^2}\;\;\;\;~~~~~\mbox{and}\;\;\; ~~~~  r_a=\frac{r(1-|a|^2)}{1-|a|^2r^2},\end{eqnarray*} respectively.

 \smallskip

\begin{lem}\cite{RK}
Let $r\in (0,1)$, $a\in \IB$ and $y\in E(x,r)$. Then, $$1-|x|^2\thickapprox
1-|y|^2\thickapprox [x,y],~~~~[a,x]\thickapprox[a,y]$$ and $$|E(x,r)|\thickapprox
(1-|x|^2)^n,$$ where $|E|$ denotes the Euclidean volume of $E$.
\end{lem}  \smallskip

The following integral estimate  will be used in the
sequel.

\begin{lem}\cite{RK1}
Suppose that $t >-1$, $\beta \in \IR$ and $x\in \IB$.  Let

$$I_{t,\beta}(x)=\int_{\IB}\frac{(1-|y|^2)^t}{[x,y]^{n+t+\beta}}dv(y).
$$
Then, $$I_{t,\beta}(x)\approx\begin{cases}
\displaystyle (1-|x|^2)^{-\beta},\;
\;\;\beta>0,\\
\displaystyle  \log\frac{1}{1-|x|^2},\;
\;\;\beta=0,\\
\displaystyle  1,\;\;\;\;\;\;\;\;\;\;\;\;\;\;\;\;\;\;\;\beta<0.\end{cases}$$
\end{lem}

Let $u\in h(\IB)$. We use the notations $$u_i(x)=\frac{\partial u}{\partial x_i}(x), \;\;\;\;\  u_{ij}(x)=\frac{\partial^2 u}{\partial x_i\partial x_j}(x),\;\;\;\;\  i,j\in \{1,2...,n\}.$$
The radial derivative $u_r$ and the second order radial derivative $u_{rr}$ of $u$ are  given by $$u_r(x)=\frac{1}{|x|}\langle\nabla u(x), x \rangle\;\;\ \mbox{and}\;\;\ u_{rr}(x)=\frac{1}{|x|^2}\sum^n_{i,j=1}x_ix_ju_{ij}(x) ,$$ respectively.
Obviously,  both $|x|u_r(x)$ and $|x|^2u_{rr}(x)$ are harmonic on $\IB$ (cf. \cite{Pa}).

We recall the following useful  inequalities concerning harmonic functions. See, for example, \cite{ABR, RK1} for the proofs.

\begin{lem}
Let $0<p < \infty$ and $u \in h(\IB)$. Then
\begin{flalign*}
&(1) ~~|u(x)|^p\lesssim\frac{1}{\varepsilon^n}\int_{\IB(x,\varepsilon)}
|u(y)|^pdv(y);&\\
&(2) ~~|\mathcal{R}u(x)|^p \lesssim\frac{1}{\varepsilon^{n+p}}  \int_{\IB(x,\varepsilon)}
|u(y)|^pdv(y);&\\
&(3) ~~|\nabla u(x)|^p \lesssim \frac{1}{\varepsilon^{n+p}}  \int_{\IB(x,\varepsilon)}
|u(y)|^pdv(y),
\end{flalign*}
where $\IB(x,\varepsilon) \subset \IB$ and $\mathcal{R}u(x)=|x|u_r(x)$.
\end{lem}  \smallskip

\section{Proofs of Theorems 1.1 and 1.2}
In this section, we shall prove  Theorem 1.1 and Theorem 1.2. In order to prove the results, we need some preparations.

\begin{lem}
Let $\beta \in (0,2)$ and $u\in h(\overline{\IB})$ with the boundary function $u^*=u|_{\mathbb{S}}$. If there exists a constant $C$ such that
\begin{equation}|u^*(b')+u^*(c')-2u^*(a')|\leq C|a'-b'|^\beta, ~ a', b', c'\in \mathbb{S},  c'=\sigma_{a'}(b'),\end{equation} then
\begin{equation} |u_{rr}(x)|\lesssim(1-|x|)^{\beta-2}, ~x\in \IB \backslash \{0\}.\end{equation}
\end{lem}  \smallskip

\bpf Let $x\in \IB \backslash \{0\}$. By (1.4),  $u$ can be written as
\begin{equation*}\label{equ2.2}
u(x)=\int_{\mathbb{S}}P(x,\xi)u^*(\xi)~d\sigma(\xi).
\end{equation*}
Differentiating $P(x,\xi)$ twice with respect to $x_i$ and $x_j$ gives
\begin{eqnarray*}\frac{\partial^2}{\partial x_i \partial x_j }P(x,\xi)&=&\frac{-2\delta_{ij}}{|x-\xi|^n}+\frac{2nx_i(x_j-\xi_j)+2nx_j(x_i-\xi_i)}{|x-\xi|^{n+2}}\\&&-\frac{n(1-|x|^2)\delta_{ij}}{|x-\xi|^{n+2}}+
\frac{n(n+2)(1-|x|^2)(x_i-\xi_i)(x_j-\xi_j)}{|x-\xi|^{n+4}},\end{eqnarray*}where $\delta_{ij}$ is the Kronecker delta.  By elementary calculations, we have
\begin{eqnarray*}|x|^2P_{rr}(x,\xi)&=&\frac{-2|x|^2}{|x-\xi|^n}+\frac{4n|x|^2(|x|^2-\langle x,\xi\rangle)}{|x-\xi|^{n+2}}-\frac{n|x|^2(1-|x|^2)}{|x-\xi|^{n+2}}\\&&+
\frac{n(n+2)(1-|x|^2)(|x|^2-\langle x,\xi\rangle)^2}{|x-\xi|^{n+4}}.\end{eqnarray*}

Set $\zeta=\sigma_{x'}(\xi)$. It follows from (1.5) that
$P_{rr}(x,\xi)=P_{rr}(x,\zeta)$, and thus
\begin{eqnarray*}u_{rr}(x)&=&\int_{\mathbb{S}}P_{rr}(x,\xi)u^*(\xi)~d\sigma(\xi)
\\&=&\int_{\mathbb{S}}P_{rr}(x,\xi)u^*(\zeta)~d\sigma(\xi)\\&=&
\frac{1}{2}\int_{\mathbb{S}}P_{rr}(x,\xi)\big(u^*(\xi)+u^*(\zeta)-2u^*(x')\big)~d\sigma(\xi) .\end{eqnarray*}
Since $\big||x|^2-\langle x,\xi\rangle\big|\leq |x||x-\xi|,$ we obtain
\begin{eqnarray*}\big|P_{rr}(x,\xi)\big|&\lesssim  &\frac{1}{|x-\xi|^{n+1}}+\frac{1-|x|}{|x-\xi|^{n+2}}.\end{eqnarray*}
This, together with the assumption (3.1), yields
\begin{eqnarray*}|u_{rr}(x)|&\lesssim& \int_{\mathbb{S}}\bigg (\frac{|\xi-x'|^\beta}{|x-\xi|^{n+1}}+\frac{(1-|x|)|\xi-x'|^\beta}{|x-\xi|^{n+2}} \bigg)d\sigma(\xi) .\end{eqnarray*}

In order to estimate the above integral, we split $\mathbb{S}$ into  the following two subsets
$$E=\{\xi\in \mathbb{S}: |\xi-x'|\leq 1-|x|\}\;\;\mbox{and}\;\; ~F=\{\xi\in \mathbb{S}: |\xi-x'|> 1-|x|\}.$$
Since $|\xi-x|\geq 1-|x|$ for all $\xi\in \mathbb{S}$, we obtain
\begin{eqnarray*} \int_{E}\bigg(\frac{|\xi-x'|^\beta}{|x-\xi|^{n+1}}+\frac{(1-|x|)|\xi-x'|^\beta}{|x-\xi|^{n+2}}\bigg )d\sigma(\xi)  &\leq& 2(1-|x|)^{-(n+1)} \int^{1-|x|}_0 \rho^{n-2+\beta}d\rho\\&\lesssim&(1-|x|)^{\beta-2}.\end{eqnarray*}
On the other hand, $|x'-\xi|\lesssim |x-\xi|$ for all $\xi \in F$, so
\begin{eqnarray*} \int_{F}\bigg(\frac{|\xi-x'|^\beta}{|x-\xi|^{n+1}}+\frac{(1-|x|)|\xi-x'|^\beta}{|x-\xi|^{n+2}}\bigg)d\sigma(\xi)  &\lesssim&  \int_{F}\frac{1}{|x-\xi|^{n+1-\beta}}d\sigma(\xi)\\&\lesssim& \int^2_{1-|x|}\rho^{\beta-3}d\rho\\&\lesssim& (1-|x|)^{\beta-2} .\end{eqnarray*}
Hence
\begin{eqnarray*}|u_{rr}(x)|&\lesssim& (1-|x|)^{\beta-2}.\end{eqnarray*}
This completes the proof.   \epf

For $u\in h(\IB)$, it is a familiar fact that the radial derivative $u_r$ controls the gradient $\nabla u$, see, e.g. \cite[Lemma 6]{Pa}. In the following lemma, we extend this result to the second order setting.

\begin{lem}
Let $\beta \in (0,2)$ and $u\in h(\overline{\IB})$. If  $u$ satisfies (3.2), then
\begin{eqnarray}|\nabla^2 u(x)|\lesssim (1-|x|)^{\beta-2}, ~ x\in \IB.\end{eqnarray}
\end{lem}
\bpf    Let $M=\max\{M_1,M_2,M_3\}$, where $M_1= \sup_{x\in \overline{\IB}}|u(x)|$, $M_2=\sup_{x\in \overline{\IB_{1/2}}}|\nabla u(x)|$, and $M_3=\sup_{x\in \overline{\IB_{5/8}}}|\nabla^2u(x)|$. Then it is easy to see that $M$ is bounded.
The remaining proof consists of the following three steps.
\smallskip

Step 1.  $|u_{r}(x)| \lesssim (1-|x|)^{\beta-2}, x\in \IB \backslash \{0\}$. \smallskip

Let $x\in \IB \backslash \{0\}$. By the mean value theorem and (3.2),
$$|u_r(x)-u_r(x/2)|\leq|x|/2\sup_{t\in [1/2,1]}|u_{rr}(tx)|\leq C (1-|x|)^{\beta-2},$$
and thus
\begin{eqnarray*}|u_{r}(x)|\leq |u_r(x/2)|+C (1-|x|)^{\beta-2} &\leq&  M_22^{2-\beta} (1-|x|)^{\beta-2}+C (1-|x|)^{\beta-2}\\&\lesssim &  (1-|x|)^{\beta-2}.\end{eqnarray*}

  Step 2.  $|(u_{ij})_{rr}(x)|\lesssim (1-|x|)^{\beta-4}, ~ |x|>1/8$.
\smallskip

 Note that $u(x), |x|u_r(x)$ and $|x|^2u_{rr}(x)$ are all harmonic on $\IB$. Let $x\in \IB$ with $|x|>\frac{1}{8}$. Setting $\varepsilon= (1-|x|)/7$ and applying Lemma 2.3 to each of these quantities gives
\begin{eqnarray}&& |u_{ij}(x)| + |(|x|u_r(x))_{ij}|+|(|x|^2u_{rr}(x))_{ij}|\notag \\ & \lesssim& \frac{1}{\varepsilon^{n+2}}\int_{\IB(x,\varepsilon)}
\big(|u(y)|+ ||y|u_r(y)|+||y|^2u_{rr}(y)|\big)dv(y)\notag \\ & \lesssim& \frac{1}{\varepsilon^{n+2}}\int_{\IB(x,\varepsilon)}
\big(|u(y)|+ |u_r(y)|+|u_{rr}(y)|\big)dv(y)\\ & \lesssim& \frac{1}{\varepsilon^{n+2}}\int_{\IB(x,\varepsilon)}
(2M+C) (1-|y|)^{\beta-2}dv(y) \notag  \\ & \lesssim&
(1-|x|)^{\beta-4}.\notag\end{eqnarray}

By direct computations, we have
\begin{eqnarray*}(|x|u_r(x))_{ij}=2u_{ij}(x)+|x|(u_{ij})_r(x),
\end{eqnarray*}
and
\begin{eqnarray*}\big(|x|^2(u_{rr}(x))\big)_{ij}=2u_{ij}(x)+4|x|(u_{ij})_r(x)+|x|^2(u_{ij})_{rr}(x).
\end{eqnarray*}
Thus,
\begin{eqnarray*}|x|^2(u_{ij})_{rr}(x)= (|x|^2u_{rr}(x))_{ij} +6u_{ij}(x)-4(|x|u_r(x))_{ij}.
\end{eqnarray*}
This, together with (3.4), yields
\begin{eqnarray*}|(u_{ij})_{rr}(x)|\lesssim (1-|x|)^{\beta-4}, ~ |x|>1/8.
\end{eqnarray*}
\smallskip
Step 3. $|u_{ij}(x)|\lesssim (1-|x|)^{\beta-2}$.

We only need to consider the case $|x|>1/2$, since  for $|x|\leq 1/2$ the assertion is trivial.  Set $x=r x'$, $r=|x|$. By the fundamental theorem of calculus,
\begin{eqnarray*}&&u_{ij}(x)+u_{ij}(x/4)-2u_{ij}(5x/8) \\&=& u_{ij}(rx')+u_{ij}(rx'/4)-2u_{ij}(5rx'/8)\\&=& \int^{3r/8}_0(3r/8-t)\Big( \frac{\partial^2 u_{ij}}{\partial r^2}((5r/8-t)x')+ \frac{\partial^2 u_{ij}}{\partial r^2}((5r/8+t)x') \Big) dt.\end{eqnarray*}
Consequently,
\begin{eqnarray*}&&\big|u_{ij}(x)+u_{ij}(x/4)-2u_{ij}(5x/8)\big|\\&\leq& \int^{3r/8}_0(3r/8-t)\Big(\Big| \frac{\partial^2 u_{ij}}{\partial r^2}((5r/8-t)x')\Big|+ \Big|\frac{\partial^2 u_{ij}}{\partial r^2}((5r/8+t)x')\Big| \Big)dt
\\&\lesssim& \int^{3r/8}_0(3r/8-t) \big(1-5r/8-t\big)^{\beta-4}dt
\\ &\lesssim&
\int^{3r/8}_0 \big( 1-5r/8-t\big)^{\beta-3} dt\\ &\lesssim &  (1-r)^{\beta-2}
.\end{eqnarray*}
Therefore,
\begin{eqnarray*}|u_{ij}(x)| &\leq & |u_{ij}(x/4)|+2|u_{ij}(5x/8)|+C (1-|x|)^{\beta-2} \\&\leq & (3M+C) (1-|x|)^{\beta-2}. \end{eqnarray*}
The proof is complete.   \epf

\begin{lem}
Let $\beta \in (0,2)$, $0<r<1$ and $u\in h(\overline{\IB})$. If $u$ satisfies (3.3), then
$$|u(\xi)+u(r\xi)-2u\big((1+r)\xi/2\big)|\lesssim (1-r)^\beta, ~ \xi\in \mathbb{S} .$$
\end{lem}
\bpf By the fundamental theorem of calculus, we have
\begin{eqnarray*}&&\big|u(\xi)+u(r\xi)-2u\big((1+r)\xi/2\big)\big|\\&\leq& (1-r)^2/4\int^1_0(1-t)\Big \{\big|\nabla^2 u\big((1+r)\xi/2+(1-r)t\xi/2\big)\big|\\&&+ \big(\big|\nabla^2 u\big((1+r)\xi/2-(1-r)t\xi/2\big)\big|\Big \}dt.
\end{eqnarray*}
Since $1-|(1+r)\xi/2+(1-r)t\xi/2|=(1-t)(1-r)/2\leq 1-|(1+r)\xi/2-(1-r)t\xi/2|$, we obtain
\begin{eqnarray*}|u(\xi)+u(r\xi)-2u\big((1+r)\xi/2\big)|
\lesssim (1-r)^\beta\int^1_0(1-t)^{\beta-1}dt \lesssim (1-r)^\beta,
\end{eqnarray*}
which completes the proof.   \epf

\begin{lem}
Let $\beta \in (0,2)$, $u\in h(\overline{\IB})$ with the boundary function  $u^*=u|_{\mathbb{S}}$. If $u$ satisfies (3.3), then
\begin{eqnarray}|u^*(b')+u^*(c')-2u^*(a')|\lesssim |a'-b'|^\beta, ~ a', b', c'\in \mathbb{S},  c'=\sigma_{a'}(b').\end{eqnarray}
\end{lem}

\bpf  Without loss of generality, we may assume
that $a'=e_1, b'=(\cos\theta, \sin\theta,0,...,0), c'=(\cos\theta, -\sin\theta,0,...,0)$ and $\theta$ is a  small positive number.

Choose $l_1$ satisfying $7/8<l_1<1$ and $1-l_1\thickapprox l=|a'-b'|=2\sin(\theta/2)$. Consider the following six points $\{a_j, b_j,c_j\}, j=1,2$ in $\IB$, where
$$a_1= l_1a', \;\;~~b_1=l_1b',\;\; ~~~~c_1= l_1c',$$ $$a_2= (2l_1-1)a' ,\;\;\;b_2=(2l_1-1)b',\;\;\ c_2= (2l_1-1)c'.$$  By the triangle inequality,
 \begin{eqnarray*}
 	|u(b')+u(c')-2u(a')| \leq 2I_1+ I_2+I_3+2I_4+I_5,
 \end{eqnarray*}
 where
\begin{align*}
	I_1&= |u(a')+u(a_2)-2u(a_1)|,\\
	I_2&=|u(b')+u(b_2)-2u(b_1)|, \\
	I_3&=|u(c')+u(c_2)-2u(c_1)|,  \\
	I_4&=|u(b_1)+ u(c_1)-2u(a_1)|, ~\mbox{ and }\\
	I_5&=|u(b_2)+u(c_2) -2u(a_2)| .
\end{align*}
It follows from Lemma 3.3 that
$$|I_k|\lesssim l^\beta, k=1,2,3.$$

We now estimate $I_4$. Let $\gamma(t)=(l_1\cos t, l_1\sin t,0,...,0), \; t\in [-\theta, \theta]$. Then
\begin{eqnarray*}I_4&=& |u\circ\gamma(\theta)+ u\circ\gamma(-\theta)-u\circ\gamma(0)|
\\&=& \int^\theta_0(\theta-t)\Big( \frac{d^2}{dt^2}u\circ\gamma(t)+  \frac{d^2}{dt^2}u\circ\gamma(-t)   \Big)dt
\\ &\leq& \int^\theta_0(\theta-t)\Big( |\nabla^2u(\gamma(t))| + |\nabla^2u(\gamma(-t))| +|\nabla u(\gamma(t))| +|\nabla u(\gamma(-t))| \Big)dt .\end{eqnarray*}
By the mean value theorem and the hypothesis on $\nabla^2 u$, we have
$$|\nabla u(x)|\lesssim (1-|x|)^{\beta-2},~x\in \IB.$$
Consequently,
\begin{eqnarray*}I_4 \lesssim \int^\theta_0(\theta-t)(1-l_1)^{\beta-2}dt
\lesssim l^\beta .\end{eqnarray*}

For the estimate of $I_5$, we consider the curve $\gamma_1(t)=\big((2l_1-1)\cos t, (2l_1-1)\sin t,0,...,0)\big), \; t\in [-\theta, \theta]$. Note that $1-(2l_1-1)=2(1-l_1)\thickapprox l$, it follows from a similar argument as in the estimate of $I_4$, one obtains $I_5\lesssim l^\beta$.
Gathering all the above estimates, (3.5) follows.   \epf

\medskip

{\bf \noindent Proof of Theorem 1.1.} By Lemmas 3.1, 3.2, and 3.4, we obtain the equivalence of $(a)$ and $(b)$. It follows from \cite[Proposition 2]{Mi} that $(a)$ is equivalent to $(c)$. The proof is finished.

\medskip

{\bf \noindent Proof of Theorem 1.2.} Without loss of generality, choose $x_0=re_n\in \IB$, and let $x_0'=e_n$. From the proof of Lemma 3.1, we have
\begin{eqnarray*}\frac{\partial^2}{\partial x_i \partial x_j }P(x_0,\xi)&=&\frac{-2\delta_{ij}}{|x_0-\xi|^n}-\frac{n(1-|x_0|^2)\delta_{ij}}{|x_0-\xi|^{n+2}}+
\frac{n(n+2)(1-|x_0|^2)\xi_i\xi_j}{|x_0-\xi|^{n+4}},\end{eqnarray*}
for $1\leq i,j <n$.

Set $\zeta=\sigma_{x_0'}(\xi)$. Then it is easy to check that
$$\frac{\partial^2}{\partial x_i \partial x_j }P(x_0,\xi)=\frac{\partial^2}{\partial x_i \partial x_j }P(x_0,\zeta)$$ and
\begin{eqnarray*}u_{ij}(x_0)&=&\int_{\mathbb{S}}\frac{\partial^2}{\partial x_i \partial x_j }P(x_0,\xi)u^*(\xi)~d\sigma(\xi)
\\&=&\int_{\mathbb{S}}\frac{\partial^2}{\partial x_i \partial x_j }P(x_0,\zeta)u^*(\zeta)~d\sigma(\xi)\\&=&
\frac{1}{2}\int_{\mathbb{S}}\frac{\partial^2}{\partial x_i \partial x_j }P(x_0,\xi)\big(u^*(\xi)+u^*(\zeta)-2u^*(x_0')\big)~d\sigma(\xi).\end{eqnarray*}
Consequently,
\begin{eqnarray*}|u_{ij}(x_0)|&\lesssim&\int_{\mathbb{S}}\bigg(\frac{|\xi- x'_0|^2}{|\xi-x_0|^n}+ \frac{(1-|x_0|)|\xi- x'_0|^2}{|\xi-x_0|^{n+2}}\bigg)~d\sigma(\xi)
.\end{eqnarray*}
By a similar argument as the proof of Lemma 3.1, we obtain
\begin{eqnarray*}|u_{ij}(x_0)|&\lesssim&
 \int_{E}\bigg(\frac{|\xi- x'_0|^2}{|\xi-x_0|^n}+ \frac{(1-|x_0|)|\xi- x'_0|^2}{|\xi-x_0|^{n+2}}\bigg)d\sigma(\xi)\\&&+\int_{F}\bigg(\frac{|\xi- x'_0|^2}{|\xi-x_0|^n}+ \frac{(1-|x_0|)|\xi- x'_0|^2}{|\xi-x_0|^{n+2}}\bigg)d\sigma(\xi)\\&\lesssim& (1-|x_0|)^{-n} \int^{1-|x_0|}_0 \rho^{n}d\rho +(1-|x_0|)^{-n-1} \int^{1-|x_0|}_0
 \rho^{n}d\rho\\&&+
 \int^2_{1-|x_0|} \big(1+ (1-|x_0|)\rho^{-2}\big) d\rho\\&\lesssim& C,
\end{eqnarray*}
and thus
$$|\nabla^2(u|_T)(x_0)|\leq C(n).$$
Due to rotational symmetry, the same estimate holds for every
derivative in any tangential direction. This establishes the boundedness of $|\nabla^2(u|_T)|$.

Combining this with the assumption (1.6) gives that $|\nabla^2 u|\leq LC(n)$.
By the fundamental theorem of calculus,
     \begin{eqnarray*}|u(x+y)+u(x-y)-2u(x)|&\leq& |y|^2\int^1_0(1-t)\big(|\nabla^2 u(x+ty)|+|\nabla^2 u(x-ty)|\big)dt\\&\leq& LC(n)|y|^2, \end{eqnarray*}for all  $x\in \IB$ and for all $y$ with $x\pm y\in \IB$.
The proof is complete.  \hfill $\Box$

\section{Proofs of Theorems 1.3--1.5}

{\bf \noindent Proof of Theorem 1.3.} We first assume that (1.8) holds.  For  $x\in \IB$ and $y$ with $x\pm y\in \IB$, let $$f(y)=u(x+y)+u(x-y)-2u(x).$$
Then $f$ is harmonic on $\IB(0, 1-|x|)$. Note that $\nabla^2f(0)= 2 \nabla^2u(x)$, so by taking $\varepsilon=(1-|x|)/2$, it follows from Lemma 2.3 that
\begin{eqnarray*}
|\nabla^2u(x)|\lesssim |\nabla^2f(0)| &\lesssim&
 \frac{1}{\varepsilon^{n+2}}\int_{\IB(0,\varepsilon)}|f(y)|\,dv(y).
 \end{eqnarray*}
Note that for $y\in \IB(0,\varepsilon)$, $x\pm y\in \IB(x,\varepsilon)$. This gives
\begin{eqnarray*}(1-|x|^2)|\nabla^2u(x)| &\lesssim&
\frac{1}{\varepsilon^n}\int_{\IB(0,\varepsilon)}(1-|x+y|^2)^{1/2}(1-|x-y|^2)^{1/2}\frac{Z_u(x,y)}{|y|^2}dv(y)\\&\lesssim& I,\end{eqnarray*}
which implies $u\in \mathcal{HZ}$, where $I$ is defined as in the statement of the theorem. \smallskip

Conversely, we assume that $u\in \mathcal{HZ}$ and $x \in \IB$, and $y\neq 0$ with $x\pm y\in \IB$. By the fundamental theorem of calculus,
\begin{eqnarray*}
Z_u(x,y) &\leq&
|y|^2\int^1_{0} (1-t) \big( |\nabla^2 u(x+ty)|+ |\nabla^2 u(x-ty)|\big)\,dt
\\&\lesssim & |y|^2\int^1_{0} \left(\frac{1}{1-|x+ty|^2}+\frac{1}{1-|x-ty|^2} \right)dt\\&\lesssim & |y|^2I_0,
\end{eqnarray*}where
$$I_0= \int^1_{-1}  \frac{1}{1-|x+ty|^2}\;dt$$
Recall that $x \in \IB$ and $y$ with $x\pm y\in \IB$. So if we write $\langle x,y\rangle=|x|\cdot|y|\cos\theta$, where $\theta$ is the angle between $x$ and $y$. Then by a straightforward calculation,
$$|x|<1, ~~|x|^2+|y|^2\pm2|x|\cdot|y|\cos\theta< 1.$$
Now $I_0$ can be written as
\begin{eqnarray*}
I_0 &=&\int^1_{-1} \frac{1}{1-|x|^2-2\langle x,y\rangle t-|y|^2t^2}\; dt \\&=&\int^1_{-1} \frac{1}{1-|x|^2+\frac{\langle x,y\rangle^2}{|y|^2}-\Big(|y|t+\frac{\langle x,y\rangle}{|y|}\Big)^2}\; dt\\&=&\int^1_{-1} \frac{1}{1-|x|^2\sin^2\theta -(|y|t+|x|\cos\theta)^2}\; dt,
\end{eqnarray*}
so that letting $|y|t+|x|\cos\theta=\big(\sqrt{1-|x|^2\sin^2\theta }\big)s$ and making the change of variables, we have
\begin{eqnarray*}
I_0 &=&\frac{1}{|y|\sqrt{1-|x|^2\sin^2\theta}}I_1,
\end{eqnarray*} where
\begin{eqnarray*}
I_1&=&\int^{\frac{|x|\cos\theta+|y|}{\sqrt{1-|x|^2\sin^2\theta}}}_{\frac{|x|\cos\theta-|y|}{\sqrt{1-|x|^2\sin^2\theta}}}\frac{ds}{1-s^2}  \\&=& \frac{1}{2}\ln\frac{(\sqrt{1-|x|^2\sin^2\theta}+|x|\cos\theta+|y|)(\sqrt{1-|x|^2\sin^2\theta}-|x|\cos\theta+|y|)}{(\sqrt{1-|x|^2\sin^2\theta}-|x|\cos\theta-|y|)(\sqrt{1-|x|^2\sin^2\theta}+|x|\cos\theta-|y|)}
\\&=& \frac{1}{2}\ln\frac{(\sqrt{1-|x|^2\sin^2\theta}+|y|)^2-|x|^2\cos^2\theta}{(\sqrt{1-|x|^2\sin^2\theta}-|y|)^2-|x|^2\cos^2\theta} \\&=&
\frac{1}{2}\ln\frac{1-|x|^2+|y|^2+2|y|\sqrt{1-|x|^2\sin^2\theta}}{1-|x|^2+|y|^2-2|y|\sqrt{1-|x|^2\sin^2\theta}}
\\&=& \frac{1}{2}\ln\Big(\frac{\lambda+2}{\lambda-2}\Big),\end{eqnarray*}and $$\lambda=\frac{1-|x|^2+|y|^2}{|y|\sqrt{1-|x|^2\sin^2\theta}}.$$
Furthermore,
\begin{eqnarray*}(1-|x+y|^2)^{1/2}(1-|x-y|^2)^{1/2}&=&
\sqrt{(1-|x|^2-|y|^2)^2-4|x|^2|y|^2\cos^2\theta}\\&=&\sqrt{(1-|x|^2+|y|^2)^2-4|y|^2(1-|x|^2\sin^2\theta)}\\&=&|y|\sqrt{1-|x|^2\sin^2\theta}\sqrt{\lambda^2-4}
  .\end{eqnarray*}
Consequently,
\begin{eqnarray*}I_2&=&(1-|x+y|^2)^{1/2}(1-|x-y|^2)^{1/2}I_0=\varphi(\lambda),
\end{eqnarray*} where $$\varphi(\lambda)=\frac{1}{2}\sqrt{\lambda^2-4}\ln\Big(\frac{\lambda+2}{\lambda-2}\Big).$$
Note that $\lambda>2$.
By elementary calculations, we deduce that $$\varphi'(\lambda)=\frac{\lambda}{2\sqrt{\lambda^2-4}}\Big[\ln\Big(\frac{\lambda+2}{\lambda-2}\Big)-\frac{4}{\lambda}\Big]\geq 0~~~ \mbox{and} ~ \lim_{\lambda\rightarrow +\infty}\varphi(\lambda)=2,$$ which imply that $\varphi(\lambda)$ is non-decreasing and bounded. This proves the boundedness of $I$ defined in Theorem 1.3. \hfill $\Box$

\medskip

{\bf \noindent Proof of Theorem 1.4.} Assume that (1.9) holds.  By  the same reasoning as in the proof of the sufficiency in Theorem 1.3, we see that $u\in \mathcal{HZ}^\alpha$.

For the converse,  assume that $u\in \mathcal{HZ}^\alpha$ and $x\pm y \in E(x,r)$.  Then, it follows from Lemma 2.1 that

\begin{eqnarray*}
&&(1-|x+y|^2)^{\beta}(1-|x-y|^2)^{\alpha-\beta}Z_u(x,y)\\&\lesssim&
|y|^2(1-|x|^2)^{\alpha}\int^1_{0} (1-t) \big( |\nabla^2 u(x+ty)|+ |\nabla^2 u(x-ty)|\big)\,dt\\&\lesssim&
|y|^2\int^1_{0} (1-t)(1-|x|^2)^{\alpha} \left(\frac{1}{(1-|x+ty|^2)^\alpha}+\frac{1}{(1-|x-ty|^2)^\alpha} \right)dt\\&\lesssim&|y|^2\int^1_{0}(1-t)\;dt\\&\lesssim&|y|^2,
\end{eqnarray*} which implies the boundedness of $J$.  This completes the proof of Theorem 1.4.  \hfill $\Box$

\medskip

\subsection*{Proof of Theorem 1.5} Assume that (1.10) holds. By Lemma 2.3, for each fixed $a\in \IB$, we have

\begin{eqnarray*}
|\nabla^2 u(a)|^p \lesssim\frac{1}{(1-|a|^2)^n}\int_{E(a,1/2)}
|\nabla^2 u(x)|^pdv(x),\end{eqnarray*}
which in turn gives
\begin{eqnarray*}
(1-|a|^2)^{p\alpha}|\nabla^2 u(a)|^p \lesssim (1-|a|^2)^{p\alpha-n} \int_{E(a,1/2)}
|\nabla^2 u(x)|^pdv(x).\end{eqnarray*}
According to Lemma 2.1, we obtain
\begin{eqnarray*}
(1-|a|^2)^{p\alpha}|\nabla^2 u(a)|^p &\lesssim&  \int_{E(a,1/2)}
|\nabla^2 u(x)|^p(1-|x|^2)^{p\alpha-n}dv(x)\\&\lesssim&\int_{E(a,1/2)}
|\nabla^2 u(x)|^p(1-|x|^2)^{p\alpha-n}(1-|\varphi_a(x)|^2)^sdv(x)\\&\lesssim&\int_{\IB}
|\nabla^2 u(x)|^p(1-|x|^2)^{p\alpha-n}(1-|\varphi_a(x)|^2)^sdv(x),\end{eqnarray*}
which implies $u\in \mathcal{HZ}^\alpha$.

Conversely, assume that $u\in \mathcal{HZ}^\alpha$. For $a\in \IB$ and $0<p<\infty$, we have
\begin{eqnarray*}
&&\int_{\IB}|\nabla^2u(x)|^p(1-|x|^2)^{p\alpha-n}(1-|\varphi_a(x)|^2)^sdv(x)  \\&\lesssim& \Big \{\sup_{x\in \IB} |\nabla^2 u(x)|^p(1-|x|^2)^{p\alpha} \Big\}  \int_{\IB}
( 1-|\varphi_a(x)|^2)^s (1-|x|^2)^{-n}dv(x)\\&\lesssim&  \|u\|_{\mathcal{HZ}^\alpha} \int_{\IB}
\frac{(1-|a|^2)^s(1-|x|^2)^{s-n}}{[a,x]^{2s}}dv(x).\end{eqnarray*}
Applying Lemma 2.2 with parameters $t=s-n$ and $\beta=s$ yields
$$\int_{\IB}
\frac{(1-|a|^2)^s(1-|x|^2)^{s-n}}{[a,x]^{2s}}dv(x)\thickapprox 1.$$
Hence,  $$\int_{\IB}|\nabla^2u(x)|^p(1-|x|^2)^{p\alpha-n}(1-|\varphi_a(x)|^2)^sdv(x) \lesssim \|u\|_{\mathcal{HZ}^\alpha}.$$
The proof is complete.\hfill $\Box$

\section*{Declarations}

\subsection*{Author Contributions}
All authors contributed substantially to the paper and approved the final submitted version.

\subsection*{Conflicts of Interest}
The authors declare that there is no conflict of interest regarding the publication of this paper.

\subsection*{Data Availability Statement}
Data sharing not applicable to this article as no datasets were generated or analyzed during the current study.

\subsection*{Funding}
X. Fu was partially supported by NNSF
of China (Grant no. 12371073).  A. Rasila was partially supported by the Li~Ka~Shing Foundation STU-GTIIT Joint Research Grant (Grant no. 2024LKSFG06) and the NSF of Guangdong Province (Grant no. 2024A1515010467).

\end{document}